\documentclass{article}      

\usepackage{amssymb}
\usepackage{amsmath}
\usepackage{latexsym}
\usepackage[mathscr]{euscript}
\usepackage{graphicx}
\usepackage{amscd}
\usepackage{amsthm} 
\usepackage{wrapfig}
\newtheorem{theorem}{Theorem}
\newtheorem{corollary}[theorem]{Corollary}
  
\newtheorem{proposition}[theorem]{Proposition}

 \newtheorem{remark}[theorem]{Remark}
 
 \newtheorem{defn}[theorem]{Definition}

\newtheorem{definition}{Definition}
 
\newcommand{\bc}{\mathbb{C}}
\newcommand{\bp}{\mathbb{ P}}

\newcommand{\br}{\mathbb{R}}
\newcommand{\bq}{\mathbb{Q}}

\newcommand{\cc}{\mathcal{C}}

\newcommand{\cs}{\mathcal{S}}
\newcommand{\cg}{\mathcal{G}}
\newcommand{\cb}{\mathcal{B}}

\newcommand{\cm}{\mathcal{M}}
\newcommand{\cn}{\mathcal{N}}

\newcommand{\cdd}{\mathcal{D}}

\newcommand{\ca}{\mathcal{A}}

\newcommand{\hk}{\hookrightarrow}
\newcommand{\bg}{\bigskip}
\newcommand{\med}{\medskip}
\newcommand{\la}{\longrightarrow}
\newcommand{\bfl}{\begin{flushleft}}
\newcommand{\efl}{\end{flushleft}}

\newcommand{\bcp}{\bc \bp}
\newcommand{\msc}{\mathscr}
 \newcommand{\Kappa}{\msc K}

\newcommand{\xr}{\xrightarrow}

\newcommand{\R}{\mathbb R}

\newcommand{\Z}{\mathbb Z}

\newcommand{\C}{\mathbb C}

\newcommand\qu{/\kern-.7ex/}%

\newcommand{\fG}{\mathfrak{g}}
\newcommand{\fZ}{\mathfrak{z}}

\DeclareMathOperator{\Aut}{Aut}

\DeclareMathOperator{\Hom}{Hom}

\DeclareMathOperator{\Diff}{Diff}

\newfont{\german}{eufm10}

\begin{document}  

  \title{Universal moduli spaces of  surfaces with flat connections and cobordism theory}
 \author{Ralph L. Cohen 
  \thanks{The first author was  partially supported by a research grant from the NSF} 
 \\ Dept. of Mathematics\\ Stanford University\\  
 Stanford, CA 94305
 \and
Soren Galatius  \thanks{The second author was  partially supported by
  NSF grant DMS-0505740 and Clay Mathematics Institute} \\ Dept. of Mathematics\\ Stanford University\\  
 Stanford, CA 94305
 \and
Nitu Kitchloo \thanks{The third author was  partially supported by a research grant from the NSF} 
 \\ Dept. of Mathematics  \\
 UC San Diego \\ La Jolla, CA 92093
 }
 \date{\today}
\maketitle  
\begin{abstract} Given a semisimple, compact, connected Lie group $G$
  with complexification $G^c$, we show there is a stable range in the
  homotopy type of the universal moduli space of flat connections on a
  principal $G$-bundle on a closed Riemann surface, and equivalently,
  the universal moduli space of semistable holomorphic
  $G^c$-bundles. The stable range depends on the genus of the surface.
  We then identify the homology of this moduli space in the stable
  range in terms of the homology of an explicit infinite loop space.
  Rationally this says that the stable cohomology of this moduli space
  is generated by the Mumford-Morita-Miller $\kappa$-classes, and the
  ring of characteristic classes of principal $G$-bundles, $H^*(BG)$.
  We then identify the homotopy type of the category of one-manifolds
  and surface cobordisms, each equipped with a flat $G$-bundle. We also explain how these results may be generalized to arbitrary compact connected Lie groups. Our
  methods combine the classical techniques of Atiyah and Bott, with
  the new techniques coming out of Madsen and Weiss's proof of
  Mumford's conjecture on the stable cohomology of the moduli space of
  Riemann surfaces.  \end{abstract}

 \tableofcontents

 \section{Introduction and statement of results}

\noindent
Let $G$ be a  fixed, connected, semisimple, compact Lie group $G$.  Our goal is to study the moduli spaces  of flat $G$-connections on principal bundles over Riemann surfaces.  By allowing the complex structures of the surfaces to vary, we are able to prove a stability theorem in homology, and to study a cobordism category built out of such  moduli spaces.   

The study of  moduli  spaces of flat connections, and its connection with holomorphic bundles on Riemann surfaces goes back to the seminal work of Atiyah-Bott \cite{atiyahbott}. Given a Riemann surface $\Sigma$ without boundary, and a principal $G$-bundle $E$,  Atiyah and Bott studied  the space of holomorphic structures on the complexification $E^c$. This space admits a canonical (complex-gauge) equivariant stratification that was first described in the work of Harder-Narasimhan \cite{HN}. The open-dense stratum for this stratification is also known as the space of semistable complex structures on $E^c$. The codimension of the remaining strata grows linearly in the genus, and hence the semistable stratum approximates the whole space increasingly with genus of the curve. 

The space of holomorphic structures may be identified with the space
of principal connections on $E$. Atiyah and Bott also study the
Yang-Mills functional on this space of all connections. They show that
the Yang-Mills functional behaves like a perfect, (gauge) equivariant
Morse-Bott function, with critical subspace given by the Yang-Mills
connections. The analytical aspects of the Yang-Mills flow were not
studied by Atiyah and Bott in \cite{atiyahbott}. However, the authors
do motivate the reason why the Harder-Narsimhan stratification
represents the descending strata for the critical level sets of the
Yang-Mills functional. In particular, this suggests that the open
stratum (identified with semistable complex structures on $E^c$) must
equivariantly deform onto the space of minimal Yang-Mills connections
on $E$. This minima can be described in terms of cental connections (see the final section). In the setting where $G$ is semisimple, these Yang-Mills minima are simply the flat connections.

The Morse theoretic program suggested above was completed
by G. Daskalopoulos \cite{dask} and J.\ R\aa de
\cite{rade}. In \cite{dask} and \cite{rade} the authors succeeded in
proving the long time convergence of the Yang-Mills flow, thereby
rigorously establishing the correspondance between the Yang-Mills
moduli spaces, and the moduli spaces of semistable complex structures.

In this paper we consider the moduli space of flat connections on
bundles, parametrized over the moduli space of Riemann surfaces of a
fixed genus. We call this the {\em universal moduli space} of flat
connections.  More specifically, we let
\begin{align}\label{universal}
  \cm^G_g = \{(\Sigma, E, \omega) : \, &\Sigma \,\, \text{is a closed
    Riemann surface of genus $g$, \,
    $E \to \Sigma$}   \\
  &\text{ is a principal $G$-bundle, and} \, \, \omega \text{ is a
    flat connection on $E$}\}/\sim.  \notag
 \end{align} 
 The relation denoted by $\sim$ is induced by diagrams of the form,
  $$
 \begin{CD}
E_1  @>\tilde \phi >>   E_2 \\
 @VVV   @VVV \\
\Sigma_1 @>>\phi > \Sigma_2
 \end{CD}
 $$
 where $\tilde \phi$ is an isomorphism of $G$-bundles living over an
 orientation preserving diffeomorphism, $\phi$.  Here $\phi$ takes the
 complex strucure of $\Sigma_1$ to the complex structure of
 $\Sigma_2$, and $\tilde \phi$ pulls back the flat connection
 $\omega_2$ to the flat connection $\omega_1$.

Strictly speaking, the moduli space $\cm^G_g$ should be viewed as a
topological stack.  Alternatively, one takes the \sl homotopy orbit
space \rm of appropriate group actions defining this equivalence
relation. The details of this topology will be described in section
2.1.

Our main theorem about this moduli space is theorem~\ref{mainone}
below, which calculates the homology $H_q(\cm^G_g)$ for $2q + 4 \leq
g$.  Before stating it, let us introduce some notation.  Let $L$
denote the canonical line bundle over $\bcp^\infty$, and let
$\bcp^\infty_{-1} = (\bcp^\infty)^{-L}$ be the Thom spectrum of the virtual
inverse $-L$ (graded so that the Thom class is in $H^{-2}$).  Let
$BG_+$ denote the classifying space $BG$, with a disjoint basepoint
added.  Roughly, theorem~\ref{mainone} says that $\cm^G_g$ and the
infinite loop space $\Omega^\infty (\bcp^\infty_{-1} \wedge BG_+)$
have isomorphic homology up to degree $(g-4)/2$.  We need to be
precise about connected components, however (the two spaces have
non-isomorphic  sets of path components).  The topological type of a principal
$G$-bundle $E \to \Sigma$ is determined by the homotopy class of a map $\Sigma \to BG$.
Since $BG$ is simply connected, this in turn is determined by the
element $f_*[\Sigma] \in H_2(BG) = \pi_1(G)$.  The correspondence
$(\Sigma,E,\omega) \mapsto f_*[\Sigma]$ defines an isomorphism $\pi_0
\cm^G_g \to \pi_1(G)$.  For $\gamma \in \pi_1(G)$, we let
\begin{align*}
  \cm^G_{g,\gamma}\subseteq \cm^G_g  
\end{align*}
denote the corresponding connected component.  All components of
$\Omega^\infty(\bcp^\infty_{-1} \wedge BG_+)$ are of course homotopy
equivalent, and we let
\begin{align*}
  \Omega^\infty_\bullet(\bcp^\infty_{-1} \wedge BG_+) \subseteq
  \Omega^\infty(\bcp^\infty_{-1} \wedge BG_+)
\end{align*}
denote the component containing the basepoint.

For completeness, we describe $\pi_0 \Omega^\infty(\bcp^\infty_{-1}
\wedge BG_+)$.  The collapse maps $BG_+ \to S^0$ and $BG_+ \to BG$
define a homotopy equivalence
\begin{align*}
  \Omega^\infty(\bcp^\infty_{-1} \wedge BG_+) \xrightarrow{\simeq}
  \Omega^\infty \bcp^\infty_{-1} \times \Omega^\infty
  (\bcp^\infty_{-1} \wedge BG).
\end{align*}
It is well known (see e.g.\ \cite{madsentillmann}) that $\pi_0
\bcp^\infty_{-1} \cong \Z$.  The Hurewicz homomorphism defines an
isomorphism
\begin{align*}
  \pi_0(\bcp^\infty_{-1} \wedge BG) \cong H_2(BG) = \pi_1(G),
\end{align*}
and we have described an isomorphism $\pi_0
\Omega^\infty(\bcp^\infty_{-1} \wedge BG_+) \cong \Z \times
\pi_1(G)$.

\begin{theorem}\label{mainone} 1.  Let $G$ be a connected, compact,
  semisimple Lie group.  Then the homology group
  $H_q(\cm^G_{g,\gamma})$ is independent of $g$ and $\gamma \in
  \pi_1(G)$, so long as $2q+4 \leq g$.

  2.  For $q$ in this range,
  $$
  H_q(\cm^G_{g,\gamma}) \cong
  H_q(\Omega^\infty_\bullet(\bcp^\infty_{-1} \wedge BG_+)).
  $$
\end{theorem}

\med
\noindent \bf Remarks.  \rm 1.  The isomorphism in the theorem is induced by an
explicit map  $\cm^G_{g,\gamma} \to
\Omega^\infty(\bcp^\infty_{-1} \wedge BG_+)$, defined by a
Pontryagin-Thom construction, cf.\ \cite{madsenweiss}, \cite{GMTW},
\cite{cohenmadsen}.  The direct definition gives a map into the
component labelled $(g-1,\gamma) \in \Z \times \pi_1(G)$, rather than
the base point component.

2.  While this description of this stable homology might seem quite
complicated, the infinite loop space appearing in this theorem has
quite computible cohomology (see \cite{galatius}).  In particular the
stable rational cohomology is essentially the free, graded commutative algebra generated
by the Miller-Morita-Mumford $\kappa$-classes, and the rational
cohomology of $BG$.  See Corollary \ref{rational} below.

3.  Notice that when $G = \{id\}$, $\cm^G_g = \cm_g$ is the moduli
space of Riemann surfaces.  In this case, part one of this theorem is
the Harer-Ivanov stability theorem \cite{harer}, \cite{ivanov}.  Part
2 of this theorem is the Madsen-Weiss theorem \cite{madsenweiss},
proving the Mumford conjecture.

\bf Note. \rm The homology groups in this theorem can be taken with \sl any \rm coefficients.  Indeed the theorem is true for any  connective generalized homology theory.

 \med
 
We then go on to interpret this theorem in terms of a stability result for the homology of the  universal moduli space of semistable holomorphic bundles, and also in terms of the $\mbox{Out}(\pi_1(\Sigma_g))$-equivariant homology of the representation variety,  $\mbox{Rep} (\pi_1(\Sigma_g), G)$.  Here $\pi_1(\Sigma_g)$ is the fundamental group of a closed, connected, oriented surface $\Sigma_g$ of genus $g$, and $ \mbox{Out}(\pi_1(\Sigma_g))$ is the outer automorphism group.  See Theorems \ref{semistability} and \ref{represent} below.

\med

The second main result of the paper regards a cobordism category of
surfaces with flat connections.  We call this category $\cc^F_G$ whose
objects are closed, oriented one-manifolds $S$ equipped with
connections on the trivial principal bundle $S \times G$, and whose
morphisms are surface cobordisms $\Sigma$ between the one-manifold
boundary components, equipped with flat $G$ bundles $E \to \Sigma$
that restrict on the boundaries in the obvious way.  (See section two
for a careful definition.)  Our result is the identification of the
homotopy type of the geometric realization of this category.

\med
\begin{theorem}\label{equivalence}  There is a homotopy  equivalence,
$$
|\cc^{F}_G| \simeq \Omega^\infty(\Sigma (\bcp^\infty_{-1} \wedge BG_+)).
$$
\end{theorem}

\med In order to prove this theorem, we will compare the category
$\cc^{F}_G$ of surfaces with flat connections to the category of
surfaces with \sl any \rm connection, $\cc_G$.  That is, this category
is defined exactly as was the category $\cc^{F}_G$, except that we
omit the requirement that the connection $\omega$ on the principal
$G$-bundle $E \to \Sigma$ be flat.  We will use results of
\cite{donaldson} and \cite{MW2} to prove that the inclusion of
cobordism categories $\cc_G^F \hk \cc_G$ induces a homotopy
equivalence on their geometric realizations, and then use the results
of \cite{GMTW} to identify the resulting homotopy type.

  \med
We have chosen to work with {\em semisimple}, compact, connected Lie groups so as to ensure that the minimal Yang-Mills connections are flat. All our theorems and techniques have obvious extensions to arbitrary connected, compact Lie groups without much more effort, (see the final section of this paper).

\med
The authors are grateful to M.J. Hopkins for originally asking about the cobordism category of flat connections, and to J. Li, I. Madsen,  and R. Vakil for helpful conversations about this work.  

\section{A stability theorem for the universal moduli space of flat connections}
 
 \subsection{The main theorem}
 
 The goal of this section to prove Theorem \ref{mainone}, as stated in
 the introduction.  We begin by describing the topology of the
 universal moduli space $\cm^G_g$ more carefully.

 Recall that if $H$ is a group acting on a space $X$, the homotopy
 orbit space $X\qu H$ is defined to be the orbit space, $(EH \times
 X)/ H$, where $EH$ is a contractible space with an action of $H$,
 such that the projection $EH \to E/H$ is a principal $H$-bundle
 (equivalently: it has local sections).  Note that if the action of
 $H$ on $X$ is free (and $X \to X/H$ is a principal $H$-bundle), then
 $X \qu H$ is homotopy equivalent to the geometric orbit space, $X/H$.

 \begin{definition}\label{flat} a.  Let $E$ be a principal $G$- bundle
   over a fixed closed, smooth, oriented surface $\Sigma_g$ of genus
   $g$.  Let $J(\Sigma_g)$ be the space of (almost) complex structures
   on $\Sigma_g$, let $\ca_{F}(E)$ be the space of flat connections on
   $E$, and define
   \begin{align}
     \cn(E) = J(\Sigma_g) \times \ca_{F}(E).
    \notag
 \end{align}
 
 \med b.  Let $\mbox{Aut}(E)$ be the group of bundle automorphisms
 of $E$, living over diffeomorphisms of $\Sigma_g$.  That is, an
 element of $\mbox{Aut}(E)$ is an isomorphism of principal
 $G$-bundles,
 $$
 \begin{CD}
 E   @>\tilde \phi >>   E \\
 @VVV   @VVV \\
\Sigma_g @>>\phi > \Sigma_g
 \end{CD}
 $$
 where $\phi \in \mbox{Diff} (\Sigma_g)$, the group of orientation preserving diffeomorphisms of $\Sigma_g$. 

 c. $\mbox{Aut}(E)$ has a natural action on $\cn(E)$.  A
 diffeomorphism transforms one complex structure into another, and a
 bundle automorphism pulls back a flat connection to a new flat
 connection.  We define $\cm^G_g(E)$ to be the homotopy orbit space of
 this action:
 $$
 \cm^G_g(E) = \cn^G_g(E) \qu \mbox{Aut}^G_g(E).
 $$
 
 d.  We define the  universal moduli space of flat connections,  $\cm^G_g$ to be the disjoint union,
 $$
 \cm^G_g = \coprod_{\{E\}} \cm^G_g(E)
 $$
 where the   union is taken over isomorphism classes of principal $G$-bundles $E \to \Sigma_g$.
 \end{definition}
 
 \bf Remarks. \rm a. It is not hard to see that $\mbox{Aut}(E)$
 fits into a short exact sequence:
 $$
 1 \to \cg(E) \to \mbox{Aut}(E) \to \mbox{Diff}(\Sigma_g) \to 1
 $$
 where $\cg(E)$ is the gauge group of smooth bundle automorphisms
 $\tilde \phi : E \to E$ living over the identity of $\Sigma_g$.

 b. The space $\cm^G_g(E)$ was denoted $\cm^G_{g,\gamma}$ in the
 introduction.  Here $\gamma = f_*[\Sigma] \in H_2(BG) = \pi_1(G)$ for
 a classifying map $f: \Sigma \to BG$ of the bundle $E \to G$.

 c. Instead of defining $\cm^G_g(E)$ as the homotopy orbit space of
 $\mbox{Aut}^G_g(E)$ acting on $\cn^G_g(E)$, we could consider the
 quotient as a topological stack.

 \med Let $\ca(E)$ be the affine space of \sl all \rm connections on
 the bundle $E$ (no flatness required).  Define the configuration
 space, $\cb^G_g(E)$ to be the homotopy orbit space,
\begin{equation}\label{config}
  \cb^G_g(E) = \left(J(\Sigma_g) \times \ca (E)\right)\qu \mbox{Aut}^G_g(E).
\end{equation}

\med Including $\ca_{F}(E) \hk \ca(E)$ defines a natural map

$j_g : \cm^G_g(E) \hk \cb^G_g(E)$.  The following is a straightforward
consequence of the works \cite{atiyahbott, dask, rade} (see also
\cite{TW}, section 3.1).

\med
\begin{theorem}\label{include2} The map $ j_g : \cm^G_g(E) \hk
  \cb^G_g(E)$ is $2(g-1)r$-connected. Here $g$ is the genus of
  $\Sigma_g$, and $r$ denotes the smallest number of the form
  $\frac{1}{2} \, dim(G/Q)$, where $Q \subset G$ is any proper compact
  subgroup of maximal rank .  In particular $j_g$ induces an
  isomorphism in homotopy groups and in homology groups in dimensions
  less than $2(g-1)r$.
\end{theorem}
\begin{proof}
  As done by Atiyah-Bott \cite{atiyahbott}, the space of
  $G$-connections on $E$ can be identified with the space of
  holomorphic structures on the induced complexified bundle, $E^{c} =
  E \times_G G^{c}$.  Moreover in the Atiyah-Bott stratification of
  the space of holomorphic bundles, one has that the space of flat
  connections, $\ca_{F}(E)$, is homotopy equivalent to the stratum of
  semistable holomorphic bundles \cite{dask, rade}.  By considering
  the codimension of the next smallest stratum (in the partial order
  described in \cite{atiyahbott}) one knows that the inclusion of the
  semistable stratum into the entire space of holomorphic bundles is
  $2(g-1)r$-connected.  Translating to the setting of connections,
  this says that the the inclusion,
  $$
  j : \ca_{F}(E) \hk \ca (E)
  $$
  is $2(g-1)r$-connected. Since $ \ca (E)$ is affine, and thus
  contractible, we can conclude that the space of flat connections,
  $\ca_{F}(E)$ is $2(g-1)r$-connected (i.e.\ its homotopy groups
  vanish through this dimension). Therefore the product $\cn(E) =
  J(\Sigma_g) \times \ca_{F}(E)$ is $2(g-1)r$-connected.

  Now consider the following diagram of principal $\mbox{Aut}^G_g(E)$-
  fibrations,
  $$
  \begin{CD}
    J(\Sigma_g) \times \ca_{F}(E)    @>\hk >>  J(\Sigma_g) \times \ca (E)\\
    @VVV  @VVV  \\
    \cm^G_g(E)  @>j_g >>  \cb^G_g(E)  \\
    @VVV   @VVV    \\
    B\mbox{Aut}(E) @>=>> B\mbox{Aut}(E)
  \end{CD}
  $$

  The above discussion implies that the top horizontal arrow induces
  an isomorphism of homotopy groups through dimension $2(g-1)r$.
  Applying the five-lemma to the long exact sequences in homotopy
  groups induced by the two bundles, we get that the middle horizontal
  arrow, $j_g : \cm^G_g(E) \to \cb^G_g(E) $ also induces an
  isomorphism of homotopy groups in this range.
\end{proof}

\med 
 \begin{remark} The formula for $r$ may be derived easily from
   \cite[equation 10.7]{atiyahbott}. By the formula given there, the
   connectivity of the map $j_g$ is at least $2(g-1)r$, where $r$
   denotes the minimum number (over the set of all proper parabolic
   subgroups of $G^c$) of positive roots of $G$, which are not roots
   of the parabolic subgroup. This number may be rewritten as we have
   stated above. We note that if $G$ is the special unitary group
   $SU(n)$, then the largest parabolic in $SL_n(\C)$ is
   $GL_{n-1}(\C)$. Hence the number $r$ is given by $n-1$ in this
   case.
\end{remark}
This theorem states that through a range of dimensions, the universal moduli space $\cm^G_g(E)$ has the homotopy type of the classifying space of the automorphism group, $B\mbox{Aut}^G_g(E)$.   Now observe that this classifying space has the following description.

Let $EG \to BG$ be a   smooth,  universal principal $G$-bundle, so that $EG$   is contractible with a free $G$-action.     Notice that the mapping space of smooth equivariant maps,
  $C^\infty_G(E, EG)$ is also contractible, and has a free action of the group $\cg(E)$.  The action is pointwise, and clearly has slices.  Thus one has a model for the classifying space of this gauge group,
  $$
  B(\cg(E)) \simeq C^\infty_G(E, EG)/\cg(E)  \cong C^\infty (\Sigma_g, BG)_E, $$
  where $C^\infty (\Sigma_g, BG)_E$ denotes the component of the mapping space   classifying the isomorphism class of  the bundle $E$. 
  
  Let $E(\mbox{Diff}(\Sigma_g)) \to B\mbox{Diff}(\Sigma_g)$ be a
  smooth, universal principal $\mbox{Diff}(\Sigma_g)$-bundle.  A nice
  model for $E(\mbox{Diff}(\Sigma_g))$ is the space of smooth
  embeddings, $E(\mbox{Diff}(\Sigma_g)) = \mbox{Emb}(\Sigma_g,
  \br^\infty)$.  The product of the action of $\mbox{Aut}(E)$ on
  $C^\infty_G(E,EG)$ and the action (through $\mbox{Diff}(\Sigma_g)$)
  on $E(\mbox{Diff}(\Sigma_g))$ gives a free action
  on the product $E(\mbox{Diff}(\Sigma_g)) \times C^\infty_G(E, EG)$.
  The quotient of this action is a model of the classifying space
  $B\mbox{Aut}(E)$.  But notice that this quotient is given by the
  homotopy orbit space,
  \begin{equation} \label{nicemodel}
  B\mbox{Aut}(E) \simeq E(\mbox{Diff}(\Sigma_g))
  \times_{\mbox{Diff}(\Sigma_g)} C^\infty(\Sigma_g, BG)_E
  \end{equation}
  where $\mbox{Diff}(\Sigma_g)$ acts on $C^\infty(\Sigma_g, BG)_E$ by
  precomposition.  We therefore have the following corollary to
  Theorem \ref{include2}.
  
  \med
  \begin{corollary}\label{connect}  There is a natural map 
  $$
  \tilde j_g : \cm^G_g(E) \to E(\mbox{Diff}(\Sigma_g))
  \times_{\mbox{Diff}(\Sigma_g)} C^\infty(\Sigma_g, BG)_E
  $$
  which is $2(g-1)r$-connected.  By taking the disjoint union over
  isomorphism classes of $G$-bundles $E$, we then have a map
  $$
  \tilde j_g : \cm^G_g \to E(\mbox{Diff}(\Sigma_g))
  \times_{\mbox{Diff}(\Sigma_g)} C^\infty(\Sigma_g, BG)
  $$
  which is $2(g-1)r$-connected.
  
  \end{corollary}

  \med We recall from \cite{cohenmadsen} that the space
  $E(\mbox{Diff}(\Sigma_g)) \times_{\mbox{Diff}(\Sigma_g)}
  C^\infty(\Sigma_g, X)$ can be viewed as the space of smooth surfaces
  in the background space $X$ in the following sense.  As in
  \cite{cohenmadsen}, define
  \begin{align}
    \cs_g(X) = &\{(S_g, f):  \text{where $S_{g } \subset \br^\infty  $ is a smooth oriented surface of genus $g$ and } \notag  \\
    & \text{and $f : S_{g } \to X$ is a smooth map.}  \} \notag
  \end{align}
  The topology was described carefully in \cite{cohenmadsen}, which used
the embedding space $\mbox{Emb}(\Sigma_g, \br^\infty)$ for
$E(\mbox{Diff}(\Sigma_g))$.  In particular, $\cs_g(BG)$ is a model for
$E(\mbox{Diff}(\Sigma_g)) \times_{\mbox{Diff}(\Sigma_g)}
C^\infty(\Sigma_g, BG)$, and therefore corollary~\ref{connect} defines
a $2(g-1)r$-connected map $\tilde j_g : \coprod_{[E]} \cm^G_g(E)
\xr{\simeq} \cs_g(BG)$.  Again, $S_g(X)$ need not be connected.  The
correspondence $(S,f) \mapsto f_*[S]$ defines an isomorphism
$\pi_0S_g(X)\cong H_2(X) = \pi_2(X)$.  For $\gamma \in H_2(X)$, we let
\begin{align*}
  S_{g,\gamma}(X) \subseteq S_g(X)
\end{align*}
be the corresponding connected component.

Now in \cite{cohenmadsen}, the stable topology of $\cs_{g,\gamma}(X)$
was studied, for a simply connected space $X$.  The following is the
main result of \cite{cohenmadsen}.

\med
\begin{theorem}\label{cohenmadsen} For $X$ simply connected, the
  homology group $H_q(\cs_g(X))$ is independent of $g$ and $\gamma$,
  so long as $2q+4 \leq g$.  For $q$ in this range,
  $$
  H_q(\cs_{g,\gamma}(X)) \cong
  H_q(\Omega^\infty_\bullet(\bcp^\infty_{-1} \wedge X_+)).
  $$
\end{theorem}

Notice that since $G$ is assumed to be a compact, connected Lie group,
$BG$ is simply connected, so we can apply theorem~\ref{cohenmadsen}.

\med
\bf Note. \rm The homology groups in this theorem can be taken with \sl any \rm coefficients.  Indeed the theorem is true for any connective generalized homology theory.

\med
We now observe that if let $X = BG$, and put 
  Corollary \ref{connect} and Theorem \ref{cohenmadsen} together, Theorem \ref{mainone} follows. \hfill \openbox.

\med
In \cite{cohenmadsen}, the stable rational cohomology of the spaces $\cs_g(X)$ was described.  This then gives the stable rational cohomology of the universal moduli space, $\cm_g^G$.  This stable cohomology is generated by the Miller-Morita-Mumford $\kappa$-classes, and the rational cohomology of $BG$.  For the sake of completeness, we state this result more carefully, and give a geometric description of how these generating classes arise.

For a graded vector space $V$ over the rationals, let $V_+$ be
positive part of $V$, i.e.\
\begin{align*}
  V_+ = \bigoplus _{n = 1}^\infty V_n.
\end{align*}
Let $A(V_+)$ be the free graded-commutative $\bq$-algebra generated by
$V_+$ \cite{milnormoore}.  Given a basis of $V_+$, $A(V_+)$ is the
polynomial algebra generated by the even dimensional basis elements,
tensor the exterior algebra generated by the odd dimensional basis
elements.  Let $\Kappa$ be the graded vector space
$H^*(\bcp^\infty_{-1};\bq)$.  It is generated by one basis element,
$\kappa_i$, of dimension $2i$ for each $i \geq -1$.  Explicitly,
$\kappa_{-1}$ is the Thom class, and $\kappa_i =
c_1^{i+1}\kappa_{-1}$, for $c_1 = c_1(L) \in H^2(\bcp^\infty)$.
Consider the graded vector space
\begin{align*}
  V = H^*(\bcp^\infty_{-1} \wedge BG_+; \bq) = \Kappa \otimes H^*(BG;\bq).
\end{align*}
Then $H^*(\Omega^\infty_\bullet (\bcp^\infty_{-1} \wedge BG_+); \bq)$
is canonically isomorphic to $A(V_+)$, and we get the following
corollary of the stable rational cohomology $H^*(\cs_{g,\gamma}(X))$
given in \cite{cohenmadsen} and Corollary \ref{connect} above.

\med
\begin{corollary}\label{rational} There is a homomorphism of algebras,
  $$
  \Theta: A((\Kappa \otimes H^*(BG; \bq))_+) \la H^*(\cm^G_{g,\gamma}
  ; \bq)
  $$ 
  which is an isomorphism in dimensions less than or equal to
  $(g-4)/2$.
\end{corollary}

 \med
 Given an element $\alpha \in H^*(BG; \bq)$,  we describe the image
 $$
 \Theta (\kappa_i \otimes \alpha) \in H^*(\cm^G_g; \bq).
 $$
 
 Consider the universal surface bundle over $\cm_g^G$:
 $$
 \Sigma_g \to \cm^G_{g,1} \xr{p} \cm^G_g.
 $$
 Here $\cm^G_{g,1} = \{(\Sigma, E, \omega, x), \, \text{where} \, x \in \Sigma \}/\sim$.
 In other words, a point in $\cm^G_{g,1}$ is a point in the universal moduli space $\cm^G_g$, together with a marked point in $\Sigma$.  The topology of $\cm^G_{g,1}$ is defined in the obvious way, so that the projection map $p :  \cm^G_{g,1} \xr{p} \cm^G_g$ is a fiber bundle.  
  
  Notice that the space $\cm^G_{g,1}$ has two canonical bundles over it.  The first is the ``vertical tangent bundle",  $T_{vert}\cm^G_{g,1}$.  This is an oriented, two dimensional vector bundle, whose fiber over $(\Sigma, E, \omega, x)$ is the tangent space $T_x \Sigma$. 
   The second canonical bundle is a principal $G$-bundle, $E^G_{g,1} \to \cm^G_{g,1}$, whose fiber over $(\Sigma, E, \omega, x)$ is $E_x$.  
   
   View a class $\alpha \in H^*(BG; \bq)$ as a characteristic class for $G$-bundles.  Then
   $\alpha (E^G_{g,1}) \in H^*(\cm^G_{g,1}; \bq)$ is a well defined cohomology class.
   Similarly,  since $T_{vert}\cm^G_{g,1}$ is an oriented, two dimensional bundle, it has a well defined Chern class $c_1 \in  H^2(\cm^G_{g,1}; \bq)$.  One then defines
   $ 
 \Theta (\kappa_i \otimes \alpha) \in H^*(\cm^G_g; \bq)$ to be the image under integrating along the fiber,
 $$
  \Theta (\kappa_i \otimes \alpha) = \int_{fiber} c_1^{i+1} \cup \alpha (E^G_{g,1}).
  $$
  
  \med
  \noindent
  \bf Remarks. \rm   1. The smoothness of the moduli spaces, $\cm^G_g$ and $\cm^G_{g,1}$ have not been discussed, so that fiberwise integration has not been justified.  However,
  as described in \cite{cohenmadsen} and \cite{madsentillmann}, the Pontrjagin-Thom
  construction, which realizes fiberwise integration in the smooth setting, is well defined, and gives the definition of the map $\Theta$ we are using.
  
  2. When $\alpha = 1$, $\Theta (\kappa_i)$ is exactly the
  Miller-Morita-Mumford class coming from $H^*(B\mbox{Diff}(\Sigma_g);
  \bq) = H^*(\cm_g; \bq)$, the cohomology of the moduli space of
  Riemann surfaces.

3.  Notice that the above formula makes good sense, even when $i = -1$, in that
$\Theta (\kappa_{-1}\otimes \alpha) =  \int_{fiber}  \alpha (E^G_{g,1}).$

   \subsection{Applications to semistable bundles and surface group representations}
      
      We now  deduce two more  direct corollaries of Theorem \ref{mainone} that stem from the close relationship between the space of flat connections,
      the space of semistable holomorphic bundles on a Riemann surface, and the space of representations of the fundamental group of the surface.

      \med As above, let $\Sigma_g$ be a fixed closed, oriented,
      smooth surface of genus $g$, and let $J(\Sigma_g)$ be the space
      of (almost) complex structures on $\Sigma_g$.  Let $E$ be a
      principal $G$ bundle over $\Sigma_g$, where as before, $G$ is a
      compact, connected, semisimple Lie group.  For a fixed $J \in
      J(\Sigma_g)$, let $\cc_{ss}^J(E)\subset \cc^J(E)$ be the space
      of semistable $G^c$- holomorphic bundles inside the full affine
      space of all holomorphic structures on the bundle, $(E, J)
      \times_G G^c$.  Define the space
     \begin{equation}\label{hol}
     \cc_{ss}^g(E) = \{(J, B): \, J \in J(\Sigma_g), \, \text{and} \, B \in \cc^J_{ss}(E)\}. 
     \end{equation}
     Let $\cc^g(E)$ be the full space of holomorphic bundles,
     $\cc^g(E) = \{(J, B): \, J \in J(\Sigma_g), \, \text{and} \, B
     \in \cc^J (E)\}.$ As before $\mbox{Aut}^G_g(E)$ acts on
     $\cc^g(E)$ with the semistable bundles $\cc^g_{ss}(E)$ as an
     invariant subspace.  Namely, an automorphism $(\tilde \phi,
     \phi)$ of the principal bundle,
   $$
 \begin{CD}
 E   @>\tilde \phi >>   E \\
 @VVV   @VVV \\
 \Sigma_g @>>\phi > \Sigma_g
 \end{CD}
 $$
 pulls back a holomorphic structure on $(E,J) \times_G G^c$ to a
 holomorphic structure on the bundle $(E, g^*(J)) \times_G G^c$. We
 define the universal moduli space of semistable holomorphic
 structures, and the universal moduli space of all holomorphic
 structures: $$ \cm^g_{ss}(E) = \cc^g_{ss}(E)\qu \mbox{Aut}^G_g(E)
 \quad \cdd^g(E) = \cc^g (E)\qu \mbox{Aut}^G_g(E)$$ to be the homotopy
 orbit spaces.  Notice that the space of all holomorphic structures
 $\cc^g_{ss}(E)$ is contractible, so the homotopy orbit space
 $\cdd^g(E) $ is a model of the classifying space, $B
 \mbox{Aut}^G_g(E) \simeq E(\mbox{Diff}(\Sigma_g))
 \times_{\mbox{Diff}(\Sigma_g)} C^\infty(\Sigma_g, BG)_E$.
 
 \begin{remark} To be faithful to the literature, we must define the universal moduli space of semistable holomorphic structures as the homotopy orbit space with respect to the automorphism group that extends $\mbox{Diff}(\Sigma_g)$ by the {\em complexified} Gauge group. However, this will not change the homotopy type of the space.
\end{remark}

 We now have the following:
    
     \med
     \begin{theorem}\label{semistable}  The inclusion of the universal moduli space of semistable holomorphic bundles into all holomorphic bundles,
     $$
     \cm^g_{ss}(E)  \hk \cdd^g(E) \simeq E(\mbox{Diff}(\Sigma_g)) \times_{\mbox{Diff}(\Sigma_g)} C^\infty(\Sigma_g, BG)_E
     $$
     is a $2(g-1)r$-connected map.
     \end{theorem}
     
     \med
     \begin{proof}  This follows from Theorem \ref{include2} together with the equivariant homotopy equivalence between the space of flat connections and  Atiyah-Bott's semistable stratum of the space of holomorphic bundles \cite{atiyahbott, dask, rade}.  
          \end{proof}  We therefore have the following  stability theorem for the universal moduli space of semistable holomorphic bundles as a corollary to Theorem \ref{cohenmadsen}.
          
       \begin{theorem}\label{semistability}
  Let $G$ be a connected, compact Lie group. Then the homology group
  of the universal moduli space of semistable holomorphic bundles,
  $H_q(\cm^g_{ss}(E) )$ is independent of $g$ and $E$, so long as
  $2q+4 \leq g$.  For $q$ in this range,
  $$
  H_q(\cm^g_{ss}(E) ) \cong H_q(\Omega^\infty_\bullet
  (\bcp^\infty_{-1} \wedge BG_+)).
  $$
  \end{theorem}

We conclude with an application to the space of representations of the fundamental group.   Choose a fixed basepoint $x_0 \in \Sigma_g$, and let $\pi=\pi_1(\Sigma_g, x_0)$ be the fundamental group based  at that point. Let $\mbox{Hom}(\pi, G)$ denote the space of homomorphisms from $\pi$ to $G$ . We topologize this space as a subspace of the mapping space, $\mbox{Map}(\pi, G)$.  

Let $\mbox{Aut}(\pi)$ denote the the group of homotopy classes of
basepoint preserving, orientation preserving, homotopy equivalences of
$\Sigma_g$. As suggested by the notation, we may identify
$\mbox{Aut}(\pi)$ with the subgroup of automorphisms of the
fundamental group that acts by the identity on $H^2(\pi, \Z) =
\Z$. The group $\mbox{Aut}(\pi)$ acts on the space of homomorphisms
$\mbox{Hom}(\pi, G)$, by precomposition.  This action descends to an action of the \sl outer \rm automorphism group,   $\mbox{Out}(\pi) =
\mbox{Aut}(\pi)/\mbox{Inn}(\pi)$  on the strict quotient  variety, $\Hom(\pi, G)/G$,  where $G$ acts      by conjugation.  Here $\mbox{Inn}(\pi)$ is the normal subgroup of \sl inner \rm automorphisms. 
   We now study how this action of $\mbox{Out}(\pi)$ lifts to an action of the orientation preserving  diffeomorphism group, $\Diff (\Sigma_g)$  on the  \sl homotopy quotient \rm space,
 $$
\mbox{Rep}(\pi, G) = \mbox{Hom} (\pi, G) \qu G.
$$

Recall that holonomy defines a homeomorphism from the space of based
gauge equivalence classes of flat connections to the corresponding
component of the space of homomorphisms:
$$
h: \ca_{F}(E) / \cg_0(E)   \xr{\cong} \mbox{Hom} (\pi, G)_E
$$
where $\cg_0(E)$ is the based gauge group (which fixes a fiber pointwise), and acts freely on $\ca_{F}(E)$ (see \cite{atiyahbott}).  This holonomy map is  $G$-equivariant, where $G$ acts as usual on $\mbox{Hom}(\pi, G)$ by conjugation, and on the space of flat connections, $\ca_{F}(E) / \cg_0(E)$, it acts by identifying $G$  as the quotient group $G = \cg(E)/\cg_0(E)$,
and  by using the action of the full gauge group $\cg(E)$ on $\ca_{F}(E)$. 
We therefore have a  homeomorphism,
$$
\ca_{F}(E)\qu \cg(E) = EG \times_G  \left( \ca_{F}(E)/\cg_0(E) \right)  \xr{\cong} EG \times_G \mbox{Hom} (\pi, G)_E = \mbox{Rep}(\pi, G)_E,
$$
where the subscript denotes the representations that correspond to the
holonomy map for $E$.  
Notice that we can take an alternative model of $\ca_{F}(E)\qu \cg(E)$ as 
$$
\ca_{F}(E)\qu \cg(E) \simeq (\ca_{F}(E) \times E(\Aut(E)) / \cg (E)
$$
which has a residual action of $\Aut(E)/\cg (E) \cong \Diff (\Sigma_g)$.  Thus the representation space,  $\mbox{Rep}(\pi, G)_E,$ is homotopy equivalent to a space with a $\Diff (\Sigma_g)$ action, and this   action clearly lifts the action of $\mbox{Out}(\pi)$ on the honest  quotient space, $ \ca_{F}(E)/\cg(E) \cong \Hom(\pi, G)/G$.   Furthermore,  for 
 genus $g \geq 2$, this diffeomorphism
group has contractible path components, and so is homotopy equivalent
to its discrete group of path components, the mapping class group,
$\mbox{Diff}(\Sigma_g) \simeq \Gamma_g = \mbox{Out}(\pi)$.  We
therefore define the $\mbox{Out}(\pi)$-equivariant homology
$$
H^{\mbox{\mbox{Out}}(\pi)}_q(\mbox{Rep}(\pi, G)_E) = H^{\Diff (\Sigma_g)}_q( \ca_{F}(E) \times E(\Aut(E)) / \cg (E))
$$

\med

\begin{theorem}\label{represent} Let $g \geq 2$.  Then the
  $\mbox{Out}(\pi)$-equivariant homology of the representation variety
  is independent of the genus $g$, so long as $2q+4 \leq g$.  For $q$
  in this range,
  $$
  H^{\mbox{Out}(\pi)}_q(\mbox{Rep}(\pi, G)_E) \cong
  H_q(\Omega^\infty_\bullet(\bcp^\infty_{-1} \wedge BG_+)).
  $$
\end{theorem}
\begin{proof} $H^{\mbox{Out}(\pi)}_q(\mbox{Rep}(\pi, G)_E) =
  H^{\mbox{Diff}(\Sigma_g)}_q(\ca_{F}(E)\qu \cg(E) )$, but the latter
  group is equal to $$H_q(E\mbox{Diff}(\Sigma_g)
  \times_{\mbox{Diff}(\Sigma_g)}\ca_{F}(E)\qu \cg(E) ) =
  H_q(\ca_{F}(E)\qu \mbox{Aut}^G_g(E)) = H_q(\cm^G_g(E)).$$ The result
  follows by Theorem \ref{cohenmadsen}.
\end{proof}

\section{The cobordism category of surfaces with flat connections}

In this section we study the cobordism category $\cc_G^F$ of surfaces
equipped with flat $G$-bundles.  In definition \ref{flat} above, we
defined moduli spaces $\cm^G_g$ of pairs $(\Sigma,E)$ consisting of a
closed Riemann surface $\Sigma$ and a flat $G$-bundle $E$ over
$\Sigma$.  In this section we generalize to Riemann surfaces with
boundary.  These moduli spaces form morphisms in $\cc_G^F$, and gluing
along common boundaries define composition of morphisms.  Some care is
needed to make this well defined (composition must be associative).
We then identify the the homotopy type of its classifying space.

We first define the relevant moduli spaces.  Let $\Sigma$ be a compact
oriented 2-manifold (not necessarily connected, possibly with
boundary).  Let $J(\Sigma)$ be the space of (almost) complex
structures on $\Sigma$.  A principal $G$-bundle $E \to \Sigma$
restricts to a principal $G$-bundle $\partial E \to \partial \Sigma$.
For a flat connection $\omega$ on $\partial E$, let $\ca_F(E,\omega)$
denote the space of flat connections on $E$ which restrict to $\omega$
on $\partial E$.  Let $\mbox{Aut}(E; \partial)$ denote the group of
automorphisms of $E$, which restrict to the identity on a neighborhood
of $\partial E$.  Thus, $\mbox{Aut}(E)$ fits into an exact sequence
\begin{align*}
  1 \to \cg(E;\partial) \to \mbox{Aut}(E,\partial) \to
  \mbox{Diff}(\Sigma) \to 1,
\end{align*}
where $\cg(E;\partial)$ is the group of gauge transformations of $E$
which restrict to the identity near the boundary.  Let $\cm(E,\omega)$
be the homotopy orbit space
\begin{align}\label{moduli-sp-bdy}
  \cm(E,\omega) = (J(\Sigma) \times \ca_F(E,\omega)) \qu
  \mbox{Aut}(E;\partial)
\end{align}

\med

An imprecise definition of $\cc_G^F$ goes as follows.

\begin{defn}\label{defCobCat}
  An object is a triple $x = (S,E,\omega)$, where $S$ is a closed
  1-manifold, $E \to S$ is a principal $G$-bundle, and $\omega$ is a connection on $E$.  The space of morphisms from $x_0 = (S_0,
  E_0, \omega_0)$ to $x_1 = (S_1, E_1, \omega_1)$ is the disjoint
  union
  \begin{align*}
    \cc_G^F(x_0, x_1) = \coprod_{E} \cm(E,\omega),
  \end{align*}
  where the disjoint union is over all $E \to \Sigma$ with $\partial E
  = E_0 \amalg E_1$, one $E$ in each diffeomorphism class, and $\omega = \omega_0 \amalg \omega_1$.

\end{defn}

This definition is imprecise because it only defines the homotopy type
of the space of morphisms, not the underlying set (the homotopy
quotient involved in defining $\cm(E)$ involves a choice).  We must
give a precise, set-level description of the homotopy quotient, and
define an associative composition on the point set level.  We present
a way of doing this.  Recall that the definition of homotopy quotient
involves the choice of a free, contractible
$\mbox{Aut}(E;\partial)$-space $E(\mbox{Aut}(E;\partial))$.  As constructed in equation (\ref{nicemodel}), 
a convenient choice of this space is given by:
\begin{align*}
  E\mbox{Aut}(E;\partial) = \R_+ \times \mbox{Emb}(\Sigma_g,[0,1]
  \times \R^\infty) \times C^\infty_G(E,EG),
\end{align*}
where $\R_+$ denote the positive real numbers,
$\mbox{Emb}(\Sigma_g,[0,1] \times \R^\infty)$ denote the space of
embeddings, which restricts to embeddings of incoming and outgoing
boundaries $S_\nu \to \{\nu\} \times \R^\infty$, $\nu = 0,1$.
$C^\infty_G(E,EG)$ denotes the space of $G$-equivariant smooth maps.
Using this space in the definition of the homotopy quotient, we get
the following definition of the set of objects and the set of
morphisms.
\begin{defn}\label{defCobCat-precise}
  A point in the space of objects $Ob(\cc^{F}_G)$, is given by a
  triple $(\mbox{S},c, \omega)$, where $\mbox{S} \subset \R^{\infty}$
  is an embedded, closed, oriented one-manifold, $c: S \to BG$ is a
  smooth map, and $\omega$ is a principal connection on the pullback
  along $c$ of $EG \to BG$.

  A point in the space of morphisms $Mor(\cc^{F}_G)$, is given by the
  data: $(t, M, i, c, \sigma)$, where $t$ is a positive real number,
  $M \subset [0,t] \times \R^{\infty}$ is a 2-dimensional cobordism,
  $i \in J(\Sigma)$ is a complex structure, and $c: \Sigma \to BG$ is
  a smooth map.  Let $E \to \Sigma$ be the pullback along $c$ of the
  universal smooth $G$-bundle $EG \to BG$.  Finally, $\sigma$ is a
  flat connection on $E$.
\end{defn}

Explicitly, given elements $(j,\tau) \in J(\Sigma) \times
\ca_F(E;\omega)$, and $(t,\phi,b) \in E\mbox{Aut}(E;\partial)$, let
$M\subseteq [0,t]\times \R^\infty$ be obtained by stretching the first
coordinate of the image $\phi(\Sigma) \subseteq [0,1] \times
\R^\infty$, and letting $i$, $c$, and $\sigma$ be induced from $j$,
$b$, and $\tau$ by the identification $M \cong \Sigma$.  This defines
an $\mbox{Aut}(E;\partial)$-invariant map
\begin{align*}
  E\mbox{Aut}(E;\partial) \times (J(\Sigma) \times \ca_F(E;\omega))
  \longrightarrow \mbox{Mor}(\cc_G^F),
\end{align*}
(with $\cc_G^F$ defined as in definition \ref{defCobCat-precise})
which descends to an injection of $\mbox{Aut}(E;\partial)$-orbits
\begin{align}\label{moduli-to-morph}
  \cm(E,\omega) \longrightarrow \mbox{Mor}(\cc_G^F).  
\end{align}
Taking disjoint union over $\Sigma$'s and $E$'s, we get an
identification of the morphism spaces in definition \ref{defCobCat}
and those in definition \ref{defCobCat-precise}.  Moreover, it is now
clear how to define an associative composition rule: take union of
subsets $M_0 \subseteq [0,t_1] \times \R^\infty$ and $M_1 \subseteq
[t_1, t_1 + t_2] \times \R^\infty$.  (For this to be a smooth
submanifold, we should insist that all cobordism $M \subseteq [0,t]
\times \R^\infty$ are ``collared'' as in \cite{GMTW}.  Similarly, $c$
and $\omega$ should be constant in the collar direction.  We omit the
details.)

\begin{remark} \label{gen.case}

  1. If all connected components of $\partial \Sigma$ have non-empty
  boundary, the action of $\mbox{Aut}(E;\partial)$ on $J(\Sigma)
  \times \ca_F(E;\omega)$ is free, and we can replace the homotopy
  quotient in (\ref{moduli-sp-bdy}) with the strict quotient. For a fixed choice of complex structure on $\Sigma$, an explicit description of this strict moduli space, parametrized over all values of $\omega$ is given by \cite{MW2}:
\[ \mathcal{M}(E) = \{ (a,c,\omega) \in \mbox{G}^{2g} \times \mbox{G}^{b-1} \times \mathcal{A}(\partial E)\, | \, \prod [a_{2i},a_{2i-1}] = \prod \mbox{Ad}_{c_j}\mbox{Hol}(\omega_j) \} \]
where $c_1=1$, and $\mbox{Hol}(\omega_j)$ denotes the holonomy of the connection $\omega$ about the j-th boundary circle over a fixed basepoint. In addition, the composition map in our category may be interpreted as a suitable symplectic reduction:

Let $\Sigma$ be a connected Riemann surface, with non-empty boundary, obtained from a (possibly disconnected) Riemann surface $\hat{\Sigma}$ by gluing along two boundary components $\mbox{B}_{\pm} \subseteq \partial \hat{\Sigma}$. Let $E$ denote a bundle over $\Sigma$ obtained by identifying a bundle $\hat{E}$ on $\hat{\Sigma}$ along $\partial{\hat{E}}$ , then  one may identify the strict moduli space $\mathcal{M}(E)$ with the symplectic reduction of the gauge group $\mathcal{G}(\mbox{G} \times \mbox{B})$ acting on $\mathcal{M}(\hat{E})$. Here $\mbox{B}$ is the one-manifold that both components $\mbox{B}_{\pm}$ are identified with, and  $\mathcal{G}(\mbox{G} \times \mbox{B})$ is identified with the gauge group of the trivial bundle over $\mbox{B}$. The moment map for the $\mathcal{G}(\mbox{G} \times \mbox{B})$-action along which the symplectic reduction is carried out is given by  $\omega \mapsto \omega_+ - \omega_-$, where $\omega_{\pm}$ denote the restrictions of the connection $\omega \in \mathcal{M}(\hat{E})$ to the boundary components $\mbox{B}_{\pm}$.

  2. Forgetting the bundle $E$ gives a functor $\cc_G^F \to
  \cc_{\mathrm{SO}(2)}$, where $\cc_{\mathrm{SO}(2)}$ is the cobordism
  category or oriented 2-manifolds from \cite{GMTW}.  An oriented
  2-dimensional cobordism $\Sigma$ defines a morphism in
  $\cc_{\mathrm{SO}(2)}$.  The inverse image of $\Sigma$ in $\cc_G^F$
  is the moduli space
  \begin{align*}
    \ca_F(E) \qu \cg(E;\partial)
  \end{align*}
  of equivalence classes of flat connections on $E$ under the action
  of the \emph{gauge group} $\cg(E;\partial)$ of gauge transformations
  of $E$ relative to the boundary.  When $E$ is trivialized, this
  group is identified with the group of smooth maps $(\Sigma,\partial
  \Sigma) \to (G,e)$.

\end{remark}

\med The goal of this section is to identify the homotopy type of the
geometric realization of this category $|\cc^{F}_G|$.  (Strictly
speaking this is the geometric realization of the simplicial nerve of
the category, sometimes known as the \sl classifying space \rm of the
category.)  More specifically, our goal is to prove Theorem
\ref{equivalence} as stated in the introduction.

In order to prove this theorem, we will compare the category
$\cc^{F}_G$ of surfaces with flat connections to the category of
surfaces with \sl any \rm connection.  Namely, let $\cc_G$ be the
category defined exactly as was the category $\cc^{F}_G$, except that
we omit the requirement that the principal connection $\omega$ be
flat.

\med The inclusion of flat connections into all connections defines a
functor
$$
\iota : \cc^{F}_G  \hk \cc_G.
$$
We will observe that this inclusion restricts to the ``positive
boundary subcategories'' defined as follows.  Let $\cc$ represent
either of the cobordism categories, $\cc^{F}_G$ or $\cc_G$.  Let
$\cc_+$ denote the subcategory that has the same objects as $\cc$, but
the morphisms of $\cc_+$ are those morphisms of $\cc$ that involve
surfaces, each path component of which has a non-empty ``outgoing''
boundary.  (The ``outgoing boundary'' of a surface $\Sigma \subset
\br^\infty \times [0,t]$ is $\Sigma \cap (\br^\infty \times \{t\})$.)
An important step in proving Theorem \ref{equivalence} is the
following.

\med
\begin{proposition}\label{positive} The inclusion functor $\iota :
  (\cc^{F}_G)_+ \hk (\cc_G)_+$ induces a homotopy equivalence of
  geometric realizations,
  $$
  \iota : |(\cc^{F}_G)_+|  \xr{\simeq}  |(\cc_G)_+|.
  $$
\end{proposition}

\med This proposition will allow us to identify the homotopy type of
$|\cc^{F}_G|$, because we will be able to identify $ |(\cc_G)_+|$ with
the geometric realization of a cobordism category studied in
\cite{GMTW}.  In that paper, the authors identified the homotopy type
of a broad range of topological cobordism categories.  We will be
interested in a particular such category we call $\cc_2(BG)$, that is
defined as follows.

\begin{defn}\label{cobBG}
  The space of objects $Ob(\cc_2(BG))$ is given by pairs $(\mbox{S}, c
  )$, where $\mbox{S} \subset \R^{\infty}$ is an embedded, closed,
  oriented one-manifold, and $c : S \to BG$ is a continuous map.

  \med The space of morphisms $Mor(\cc_2(BG))$, is given by triples,
  $(t,\Sigma, c)$, where $t$ is a nonnegative real number, $\Sigma
  \subset [0,t] \times \R^{\infty}$ is an oriented cobordism, and $c:
  \Sigma \to BG$ is a continuous map.  The embedded surface is
  collared at the boundaries as in \cite{GMTW}. In particular,
  $\Sigma_0 = \Sigma \cap (\br^\infty \times \{0\})$ and $\Sigma_t =
  \Sigma \cap (\br^\infty \times \{t\})$ are smoothly embedded,
  oriented one-manifolds. Again, morphisms are assumed to be collared, and composition is
  defined by gluing of cobordisms and maps.  (See \cite{GMTW} for
  details.)
\end{defn}

The homotopy type of the geometric realization $|\cc_2(BG)|$ was
determined in \cite{GMTW}.  Namely, the following was proved there.

\med
\begin{theorem}(\cite{GMTW})\label{GMTW} a.  Let $(\cc_2(BG))_+$
  denote the positive boundary subcategory.  Then the inclusion
  functor, $ (\cc_2(BG))_+ \hk \cc_2(BG) $ induces a homotopy
  equivalence of geometric realizations,
  $$ 
  |(\cc_2(BG))_+| \xr{\simeq} |\cc_2(BG)|.
  $$
  
  \med b.  There is a homotopy equivalence,
  $$
  |\cc_2(BG)| \simeq \Omega^\infty(\Sigma (\bcp^\infty_{-1} \wedge
  BG_+)).
  $$\hfill\openbox
\end{theorem}

Because of this theorem, Theorem \ref{equivalence} will follow from
Proposition \ref{positive}, and the following two results.

\med
\begin{proposition}\label{cc2}
  The functor $\cc_G \to \cc_2(BG)$ which on the level of morphisms is
  given by $(t,\Sigma,j,c,\omega) \mapsto (t,\Sigma,c)$ induces
  homotopy equivalences of geometric realizations
  \begin{align*}
    |\cc_G| &\simeq |\cc_2(BG)|\\
    |(\cc_G)_+| &\simeq |(\cc_2(BG))_+|.
  \end{align*}
\end{proposition}

\med
\begin{proposition}\label{pbdry} The inclusion of the positive
  boundary subcategory,
  $(\cc^{F}_G)_+ \hk \cc^{F}_G $ induces a homotopy equivalence of
  geometric realizations,
  $$
  |(\cc^{F}_G)_+| \simeq |\cc^{F}_G|.
  $$
\end{proposition}

 \med
\subsection{The positive boundary subcategories}

\med In this subsection we prove Propositions \ref{positive} and
\ref{cc2}.  In view of Theorem \ref{GMTW}, this will imply that there
is an equivalence of the geometric realization of the positive
boundary subcategory,
\begin{equation}\label{posequivalence}
  |(\cc^{F}_G)_+| \simeq \Omega^\infty(\Sigma (\bcp^\infty_{-1} \wedge BG_+)).
\end{equation}
Thus the proof of Theorem \ref{equivalence} would then be completed
once we prove Proposition \ref{pbdry}, which we will do in the next
subsection.

\med \sl Proof of Proposition \ref{positive}.  \rm Morphisms in both
categories are given by tuples $(t,\Sigma,j,c,\omega)$, where $\omega$
is a connection on a principal $G$ bundle $E \to \Sigma$.  The only
difference between the two categories is that in one of them, $\omega$
is required to be flat.  Let $(t,\Sigma,j,c)$ be fixed.  We prove that
under the ``positive boundary'' assumption on $\Sigma$, the inclusion
of flat connections into all connections,
\begin{align*}
  \ca_F(E) \to \ca(E),
\end{align*}
is a homotopy equivalence.

The ``positive boundary'' assumption implies that no connected
component of $\Sigma$ is a closed 2-manifold.  Hence $\Sigma$
deformation retracts onto its 1-skeleton $X \subseteq \Sigma$.  Choose
a 1-parameter family $\phi_t: \Sigma \to \Sigma$ of smooth maps which
start at the identity, and such that $\phi_1$ retracts $\Sigma$ onto
its 1-skeleton.  We can lift this family to a 1-parameter family
$\Phi_t: E \to E$ of maps of principal $G$-bundles with $\Phi_0$ the
identity.  Then we can let $\omega_t$ be the connection on $E$
obtained by pullback along $\Phi_t$.  The curvature of $\omega_t$ can
be computed by naturality:
\begin{align*}
  F_{\omega_t} = (\phi_t)^*(F_\omega),
\end{align*}
and therefore $\omega_1$ is flat, because $\phi_1$ has one-dimensional
image and the curvature is a two-form.  Thus the identity map of
$\ca(E)$ is homotopic to a map into $\ca_F(E)$, and since $\ca(E)$ is
contractible, $\ca_F(E)$ is contractible too.

We have proved that the functor $(\cc_G^F)_+ \to (\cc_G)_+$ induces a
homotopy equivalence on morphisms or, in other words,
\begin{align*}
  N_1(\cc_G^F)_+ \to N_1(\cc_G)_+.
\end{align*}
For $k \geq 2$, the argument is similar: An element in $N_k(\cc_G)_+$
is given by $(\Sigma, j,c,\omega)$ as before, together with a
$k$-tuple $(t_1, \dots, t_k)$ of positive real numbers.  Here $\Sigma
\subseteq [0,t]\times \R^\infty$ with $t = t_1 + \dots + t_k$.  The
same procedure as for $k = 1$ gives a path from $\omega$ to a flat
connection.  Therefore the functor induces homotopy equivalences on
$k$-nerves for all $k$, and hence on the geometric realization.\hfill
\openbox

\begin{remark} Notice that the above proof holds for any Lie group. In particular, it shows that the space of flat connections is gauge equivariantly contractible, for any principal bundle over a connected Riemann surface with non-empty boundary. 
\end{remark}

   \bg
   
   We now go about proving Proposition \ref{cc2}.
   
   \med \sl Proof of Proposition \ref{cc2}.  \rm The map in the
   proposition is induced by the functor $(t,\Sigma,j,c,\omega)
   \mapsto (t,\Sigma,c)$ which forgets the complex structure $j$ on
   the oriented surface $\Sigma$, and forgets the connection $\omega$
   on the principal $G$-bundle $E \to \Sigma$.  Thus the functor gives
   a fibration
   \begin{align}\label{n1}
     N_1\cc_G \to N_1\cc_2(BG)
   \end{align}
   whose fiber over $(t,\Sigma,c)$ is the space
   \begin{align*}
     \cn(\Sigma,c) = J(\Sigma) \times \ca(E),
   \end{align*}
   where $J(\Sigma)$ is the space of (almost) complex structures on
   the oriented surface $\Sigma$ and $\ca(E)$ is the space of
   connections on $E$.  But both spaces are contractible, so
   (\ref{n1}) is a homotopy equivalence.  The higher levels of the
   simplicial nerve are completely similar, and we get a homotopy
   equivalence of geometric realizations.
   \hfill\openbox

 \bg
 \sl Proof of Proposition \ref{pbdry}.  \rm

 In the case where $G$ is the trivial group (in other words, omit the
 flat $G$-bundle $E \to X$ from definition~\ref{defCobCat}), we recover 
 the cobordism category $\cc_d$ from \cite{GMTW}, when $d=2$.  The positive
 boundary subcategory is denoted $\cc_d^\partial$.  In \cite[Section
 6]{GMTW}, it is proved that the inclusion $|\cc_d^\partial| \to
 |\cc_d|$ is a weak homotopy equivalence when $d \geq 2$. (In \cite{GMTW} the notation $B\cc$ was used to denote
 the geometric realization of the nerve of a category $\cc$, rather than $|\cc|$.)  The proof of
 proposition \ref{pbdry} will follow \cite[Section 6]{GMTW} very
 closely.  We first recall an outline of that argument.  Let $\cc =
 \cc_{\{e\}}$, the cobordism category in the case $G$ is the trivial
 group.

In \cite{GMTW}  a functor   $D$ from smooth manifolds to sets was defined, where
$D(X)$ is the set of smooth manifolds $W\subseteq X \times \R \times
\R^\infty$ such that the projection $(\pi,f): W \to X \times \R$ is
proper, and the projection $\pi: W \to X$ is a submersion with
2-dimensional fibers.  A \emph{concordance} is an element $W \in D(X
\times \R)$; in that case the restrictions to $X \times \{0\}$ and $X
\times \{1\}$ are called \emph{concordant}.  This is an equivalence
relation on $D(X)$, and the set of equivalence classes is denoted
$D[X]$.    The equivalence $|\cc |\to |\cc^\partial |$ was proved in  \cite{GMTW} 
by proving two natural isomorphisms:
\begin{align}
  D[X] &\cong [X, |\cc|]\label{eq:1}\\
  D[X] &\cong [X, |\cc^\partial |]\label{eq:2}
\end{align}
The first is proved as follows (again, in outline).  Given an element
$W \in D(X)$ and a point $x \in X$, let $W_x \subseteq \R \times \R^\infty$ be the
$d$-manifold $W_x = \pi^{-1}(x)$ and let $f_x: W_x \to \R$ denote the
projection to the first factor.  A choice of regular value $a \in \R$
for $f_x$ defines an object $(f_x)^{-1}(a)$ of $\cc$, and if $a_0 <
a_1$ are both regular values, then $(f_x)^{-1}([a_0,a_1])$ is a
morphism in $\cc$ between the two corresponding objects.  This is used
to define a map from left to right in~(\ref{eq:1}) which is an
isomorphism.

To construct~(\ref{eq:2}), we need to ensure that only morphisms
satisfying the positive boundary condition arise as
$f^{-1}([a_0,a_1])$.  At the heart of this is \cite[lemma 6.2]{GMTW},
which constructs a continuous family of pairs $(K_t, f_t)$, $t \in \R$
consisting of a $d$-manifold $K_t$ containing the open subset $U =
\R^d - D^d \subseteq K_t$, and a smooth function $f_t: K_t \to \R$
which is constant on $U$ and proper when restricted to $K_t - U$.
Furthermore $K_0 = \R^d$ and $f_0$ is constant, and $K_1 = \R^d -
\{0\}$ and $f_1(x)$ goes to infinity as $x \to 0$.

Now let $W \in D(X)$ and let $a_0 < a_1$ be regular values of $f_x:
W_x \to \R$.  If we are lucky, $f_x^{-1}([a_0,a_1])$ already satisfies
the positive boundary condition.  If not, let $Q \subseteq
f_x^{-1}([a_0, a_1])$ be a connected component not touching
$f_x^{-1}(a_1)$, and let $e: \R^d \to Q$ be an embedding.  Gluing in
the family $(K_t, f_t)$, we get a one-parameter family of pairs $(W_t,
f_t)$.  Repeat this procedure for each such component $Q$, and we get
a one-parameter family $(W_t, f_t)$ starting at $(W_x, f_x)$ at time
0, and ending at some other pair $(W_x', f_x')$ at time 1, for which
$(f_x')^{-1}([a_0, a_1])$ satisfies the positive boundary condition.
The rest of the proof in \cite[section 6]{GMTW} describes how this
construction and a (somewhat complicated) gluing procedure can be used
to construct the map in \eqref{eq:2}.

For the purposes of proving Proposition \ref{pbdry}, we need to construct a version of the ``standard'' one-parameter
family $(K_t, f_t)$ where $K_t$ is equipped with a flat $G$-bundle
$E_t \to K_t$, which is specified over $K_0$.  Fortunately, this is
not hard.  Namely, the proof of \cite[lemma 6.2]{GMTW} also constructs
an immersion $j_t: K_t \to K_0$ which is the identity for $t = 0$.
Then any given any flat $G$-bundle  $E_0 \to K_0$, it has a canonical extension    to a flat $G$-bundle
$E_t = (j_t)^*(E_0)$ over $K_t$.  With this extension of \cite[lemma
6.2]{GMTW} in place, the rest of the proof in \cite[section 6]{GMTW}
applies verbatim if we add flat $G$-bundles to every surface in sight.
\hfill \openbox

This completes the proof of Theorem \ref{GMTW}.

\subsection{The universal moduli space and the loop space of the
  cobordism category}

We connect our two main theorems (Theorem \ref{mainone} and Theorem
\ref{equivalence}).  First note that the empty set $\emptyset$ is
one-manifold, and thereby an object in $\cc^{F}_G$.  In the case
$\partial \Sigma = \emptyset$, the definition of $\cm(E;\omega)$ in
(\ref{moduli-sp-bdy}) agrees with our previous definition of
$\cm_g^G(E)$, and as in (\ref{moduli-to-morph}) above, we get a map
\begin{align*}
  \cm_g^G(E) \to \cc_G^F(\emptyset,\emptyset),
\end{align*}
which identifies $\cm_g^G(E)$ as the (open and closed) subset of
$\cc_G^F(\emptyset,\emptyset)$ where the topological type of $M$ is
fixed to be that of a connected genus $g$ surface, and the homotopy
class of $c: M \to BG$ is fixed to be that of a classifying map for
$E$.  In particular, an element of $\cm_g^G$ defines a loop in the
classifying space $|\cc_G^F|$ that start and end at the vertex
$\emptyset \in |\cc^{F}_G|$.

\med
\begin{corollary} The induced map to the loop space of the cobordism
  category,
  $$
  \iota_g : \cm^G_g(E) \to \Omega|\cc^{F}_G|
  $$
  induces an isomorphism in homology, $H_q(\cm^G_g(E)) \xr{\cong}
  H_q(\Omega_\bullet|\cc^{F}_G|)$ for $2q \leq 2g-4$.
\end{corollary}
\begin{proof}  Consider the following diagram:
  $$
  \begin{CD}
    \cm^G_g    @>\iota_g >>  \Omega|\cc^{F}_G|  \\
    @V\tilde j_g VV     @VVj V \\
    \cs_g(BG)  @>\iota_g >> \Omega|\cc_2(BG)|
  \end{CD}
  $$
  where $\cc_2 (BG)$ is the cobordism category of oriented surfaces in
  the background space $BG$, as described in section one.  By
  construction, this diagram commutes.
    
  We proved above that the map $j : |\cc^{F}_G| \to |\cc_2(BG)|$ is an
  equivalence.  So the right hand vertical map in this diagram is an
  equivalence.  By the work of \cite{cohenmadsen} one knows that the
  bottom horizontal map $\iota_g : \cs_g(BG) \to
  \Omega_\bullet|\cc_2(BG)|$ induces an isomorphism in homology
  through dimension $\frac{g}{2}-2$.  By Corollary \ref{connect}, the
  left hand vertical map $\tilde j_g : \cm^G_g \to \cs_g(BG)$ is
  $2(g-1)r$-connected.  The result now follows.
\end{proof}

 \section{The case of a general compact connected Lie group}

  We may extend the above constructions to an arbitrary compact
  connected Lie group $G$. This involves little more than a cosmetic change in definitions. Let $\fG$ be the Lie algebra of $G$, with center $\fZ \subseteq \fG$. Given a Riemann surface $\Sigma$ with a fixed metric, we define say that a connection $\omega$ on a bundle $E$ over $\Sigma$ is {\bf central} if the curvature of $\omega$ is a constant multiple of the volume form, with values in $\fZ$. Notice
  that this makes sense since the Lie algebra $\fZ$ generates a trivial summand in the adjoint bundle $E \times_{G} \fG$. As mentioned in the introduction, the Yang-Mills functional on a closed Riemann surface achieves a minimum on the space of central connections. The curvature for any connection in this space is independent of the connection, and is given by a topological invariant of the bundle. 
  
  Given a bundle $E$ with a fixed identification of $\partial(E)$, notice that the relative gauge group $\cg(E;\partial)$ acts on the space of central connections on $E$. Let $\cm^C(E,\omega)$ denote the universal moduli space of central connections on the bundle $E$ which restrict to $\omega$ on the boundary. We may now define the category of central connections $\cc^C_G$:
  
  \begin{defn}\label{defCobCat2}
  An object of $\cc^C_G$ is a quadruple $x = (S,\sigma,E,\omega)$, where $S$ is a closed
  1-manifold, endowed with a metric $\lambda$, $E \to S$ is a principal $G$-bundle, and $\omega$ is a
  connection on $E$.  The space of morphisms from $x_0 = (S_0,\lambda_0,E_0, \omega_0)$ to $x_1 = (S_1, \lambda_1,E_1, \omega_1)$ is the disjoint
  union
  \begin{align*}
    \cc_G^C(x_0, x_1) = \coprod_{E} \cm^C(E,\omega),
  \end{align*}
  where the disjoint union is over all $E \to \Sigma$ with $\partial E
  = E_0 \amalg E_1$, one $E$ in each diffeomorphism class, and $\omega = \omega_0 \amalg \omega_1$. This morphism space is parametrized over the moduli space of Riemann surfaces with metric (which is equivalent to the one without metrics). 
\end{defn}
 As before, this definition needs to be replaced with a working definition involving the classifying space $BG$ to make it well defined.

For a closed surface $S$, the results in \cite{atiyahbott} can be invoked again to show that the space $\cm^C(E)$ approximates the classifying space of the group of total automorphisms of $E$ with connectivity given by Theorem \ref{include2}. On the other hand if $S$ is a connected surface with non-empty boundary, then it is easy to show as before that the space $\cm^C(E,\bullet)$ (where $\omega$ is free) is homotopy equivalent to the classifying space of the total automorphism group of $E$ relative to the boundary. 

Using these facts, an easy extension of previous arguments shows that our main theorems \ref{mainone}, \ref{equivalence} and \ref{semistability} remain valid for a general compact connected Lie group $G$. 

Theorem \ref{represent} about the equivariant cohomology of the representation variety may also be interpreted in this context. Recall \cite{atiyahbott}, that any central connection on a bundle $E$ over a closed Riemann surface $S$ yields a $G$-representation of the fundamental central extension $\hat{\pi}$ of $\pi_1(S)$ by $\R$. Furthermore, this representation preserves the respective centers. 

Let $\mbox{Rep}^C(\hat{\pi},G)_E = \Hom^C(\hat{\pi},G)_E \qu G$ denote the component of the variety of such central representations. It is not hard to see that $\mbox{Rep}^C(\hat{\pi},G)_E$ admits an induced action of the group of  Symplectomorphism of $S$. For genus $g \geq 2$, this group is homotopy equivalent to $\mbox{Out}(\pi)$. Theorem \ref{represent} is now true for all compact connected Lie groups once we replace $\mbox{Rep}(\pi,G)_E$ by $\mbox{Rep}^C(\hat{\pi},G)_E$. The equivariant cohomology is to be understood with respect to the group $\mbox{Symp}(S) \simeq \Gamma_g = \mbox{Out}(\pi)$.


\begin{thebibliography}{99}
\bibitem{atiyahbott}M. Atiyah, and R. Bott \emph{The Yang-Mills equations over Riemann surfaces}, Phil. Trans. R. Soc. Lond. A  \bf 308, \rm 523-615 (1982)  
\bibitem{cohenmadsen} R.L. Cohen, and I. Madsen, \emph{Surfaces in a
    background space and the homology of mapping class groups.}
  preprint, arXiv:math.GT/0601750 (2006)
\bibitem{dask} G. Daskalopoulos, \emph{The topology of the space of
    stable bundles on a compact Riemann surface}, Journal of
  Differential Geometry \textbf{36}, (1992), No. 3, 699-746.
\bibitem{donaldson} S.K. Donaldson, \emph{Boundary value problems for
    Yang-Mills fields}, Journal of Geometry and Physics \textbf{8},
  (1992), 89-122.
\bibitem{galatius} S. Galatius, \emph{Mod $p$ homology of the stable
    mapping class group}, Topology {\textbf{43}}, (2004), 1105-1132.
\bibitem{GMTW} S. Galatius, I. Madsen, U. Tillmann, and M. Weiss,
  \emph{The homotopy type of the cobordism categoriy}, Acta Math.,
  \emph{to appear}.  Available as arXiv:math/0605249.
\bibitem{HN} G. Harder, M. S. Narasimhan, \emph{On the cohomology
    groups of moduli spaces of vector bundles over curves},
  Math. Ann. 212, (1975), 215-248.
\bibitem{harer} J.L. Harer, \emph{Stability of the homology of the
    mapping class groups of orientable surfaces}
  Ann. Math. \textbf{121}, (1985), 215-249.
\bibitem{ivanov}N.V. Ivanov, \emph{On stabilization of the homology of
    Teichmuller modular groups}, Algebra i Analyz, \textbf{V. 1
    No. 3}, (1989), 110-126; \textbf{English translation:} Leningrad
  J. of Math., \textbf{V. 1, No. 3}, (1990), 675-691.
\bibitem{madsentillmann} I. Madsen and U. Tillmann, \emph{The stable
    mapping class group and $Q(\bc \bp^\infty)$}, Invent. Math.,
  \textbf{145 (3)}, (2001), 409--544.
\bibitem{madsenweiss} I. Madsen and M. Weiss, \emph{The stable moduli
    space of Riemann surfaces: Mumford's conjecture}, Ann.\ Math.\
  \textbf{165} (2007) 843--941.
\bibitem{MW1} E. Meinrenken, C. Woodward \emph{Hamiltonian loop group
    actions and Verlinde factorization}, J. Diff. Geom. \textbf{50},
  no. 3, 417-469.  (1998) (arXiv:dg-ga/9612018)
\bibitem{MW2} E. Meinrenken, C. Woodward \emph{Cobordism for
    Hamiltonian loop group actions and flat connections on the
    punctured two-sphere}, Math. Zeit. \textbf{231}, 133-168.  (1999)
  (arXiv:dg-ga/9707018)
\bibitem{milnormoore}J.W. Milnor and J.C. Moore, \emph{On the
    structure of Hopf algebras}, Annals of Math. \textbf{81} (1965),
  211-264.
\bibitem{PS} A. Pressley and G. Segal, {\em Loop Groups}, Oxford
  Math. Monographs, Clarendon Press (1986).
\bibitem{rade} J. R\aa de, \emph{On the Yang-Mills heat equation in
    two and three dimensions}, J. Reine. angew. Math. 431 (1992),
  123-163.
\bibitem{TW}C. Teleman and C. Woodward, \emph{Parabolic bundles,
    products of conjugacy classes, and quantum cohomology},
  Ann. Inst. Fourier (Grenoble) (53) (2003), 713-748.
\end{thebibliography}
 \end{document}